\documentclass[12pt]{article}

\usepackage{epsfig,wrapfig,epic}
\usepackage{amsmath}
\usepackage{amscd}
\usepackage{wrapfig}
\usepackage{blkarray}
\usepackage{hyperref}
\usepackage{tabularx}

\usepackage{amsfonts}
\usepackage{mathrsfs}
\usepackage{stmaryrd}

\usepackage{algorithm}
\usepackage{algorithmic}
\usepackage{multirow}

\usepackage{dsfont}
\usepackage{amssymb}
\usepackage{url}

\usepackage{epstopdf}
\usepackage{subfigure}

\newtheorem{example}{Example}[section]

\newtheorem{definition}{Definition}[section]
\newtheorem{theorem}{Theorem}[section]

\newtheorem{lemma}[theorem]{Lemma}
\newtheorem{remark}[theorem]{Remark}

\def\proof{\textsc{Proof. }}
\def\foorp{\hfill$\square$}

\title{Structural index reduction algorithms for differential algebraic equations
via fixed-point iteration}

\author{
Juan Tang\thanks{Chengdu Institute of Computer Applications, Chinese Academy of Sciences (CAS);
{\em Emails:} tj-123123@163.com, tangjuan0822@gmail.com}
\and Wenyuan Wu
\thanks{Automated Reasoning and Cognition Key Lab of Chongqing, CIGIT, CAS;
{\em Email:} wuwenyuan@cigit.ac.cn.}
\and Xiaolin Qin
\thanks{Chengdu Institute of Computer Applications, Chinese Academy of Sciences;
{\em Email:} qinxl@casit.ac.cn}
\and Yong Feng
\thanks{Corresponding author. Automated Reasoning and Cognition Key Lab of Chongqing, CIGIT, CAS;
{\em Email:} yongfeng@cigit.ac.cn.}
}
\date{}

\begin{document}

\maketitle
\thispagestyle{empty}

\begin{abstract}
Motivated by Pryce's structural index reduction method for differential algebraic equations (DAEs),
we show the complexity of the fixed-point iteration algorithm and
propose a fixed-point iteration method with parameters.
It leads to a block fixed-point iteration method which can be applied to
large-scale DAEs with block upper triangular structure.
Moreover, its complexity analysis is also given in this paper.


\vspace{1mm} \noindent {\bf MSC(2010):} 34A09, 65L80, 65F50, 90C05, 90C27, 90C06.

\end{abstract}

\section{Introduction}

Differential algebraic equations (DAEs) systems arise naturally in modeling many dynamical systems, such as electric circuits, mechanical systems, and spacecraft dynamics.
Based on unified multi-domain modeling techniques e.g. Modelica \cite{Modelica2012}, computers can automatically produce thousands of DAEs.
The generated DAEs have many interesting characteristics, such as large scale, high index, block structures, which are the major motivations of our work in this paper.
It is well known that a direct numerical simulation without index reduction may not be possible or may provide a bad result
\cite{Petzold82,Gear84}. Here the index of a DAE system is a key notion in the theory for measuring the distance from the given system with a singular Jacobian to the corresponding ordinary differential equations with a nonsingular Jacobian.
Various index concepts exist in the theory
of DAEs; and the one related to the structural analysis approach is the ``\textit{structural index}", which is defined in (\ref{SIndex}).
For other indices, we refer the interested readers to \cite{Brenan96,Lamour2013,Rheinboldt84}.
High-index DAE systems usually need differentiations to reveal all the system's constraints, which are crucial to determine consistent initial
conditions. This procedure is the called ``\textit{index reduction}" of DAEs. For applications of high index DAEs, see \cite{QWFR13}. Identifying all hidden constraints on formal power series solutions in the neighborhood of a given point is a key step
to construct nonsingular Jacobian of a DAE system for numerical integration. Thus, for DAE systems, index reduction is fundamental and unavoidable.

In the previous work on DAE index reduction for general DAE systems, Campbell and Gear gave a derivative-array method to reduce DAEs in \cite{Campbell95}, which may not be applicable to large-scale nonlinear systems.
Pantelides in \cite{Pantelides1988} introduced a graph-oriented method which gives a systematic way to reduce high-index of DAEs with order one to lower index, by selectively adding differentiated forms of the equations already present in the system. In \cite{Pryce2001}, Pryce developed structural analysis method which is proved to computes the same structural index with Pantelides' method and is straightforward method for analyzing the structure of DAEs of any order.
This approach is based on solving an assignment problem, which can be formulated as an integer linear programming problem. The idea was generalized to a class of partial differential algebraic equations by Wu et al \cite{WRI08}.
Recently, Pryce et al. in \cite{Pryce2012,Nedialkov2012} generalized the
structural analysis method to the DAE systems with \emph{coarse} or \emph{fine}
block triangular forms (BTF), and showed that the difference between global offsets
of signature matrix ($\Sigma$) and local offsets of each sub-block in $\Sigma$ with \emph{fine} BTF is constant.
We focus on structural index reduction method for block triangularied systems to directly calculate the smallest offsets of the system in sequential block-wise manner, and give the complexity analysis of structural index reduction algorithms for DAEs systems without or with BTF.

The rest of this paper is organized as follows.
Section \ref{completeP} briefly reviews Pryce's structural analysis method, firstly. Then we novelly and finely prove the existence and uniqueness of smallest optimal solution of Problem \ref{LPP:Dual} and show the termination of fixed-point iteration algorithm. In addition, we also give the time complexity of the algorithm which is $O(n^3+||\textbf{c}^*||_1 \cdot n^2)$ due to Theorem \ref{PryceTime} and is not given in \cite{Pryce2001}, where $n$ is the size of the system.
Section \ref{SectBFIRM} first introduces the block triangular forms (BTF) for large scale DAE systems. Based on our fixed-point iteration method with parameter, a block fixed-point iteration algorithm is proposed to find the unique smallest dual-optimal pair of the systems with BTF, and its time complexity is $O(\sum\limits_{i=1}^{\ell}{{n_i}^3+||\textbf{c}^*_i||_1 \cdot {n_i}^2})$ by Theorem \ref{BFIRM_Time}, where $\sum\limits_{i=1}^{\ell}n_i=n$ and $\ell$ is the number of the blocks on the diagonal. It is usually much better than the cost $O(n^3+||\textbf{c}^*||_1 \cdot n^2)$ without taking the advantages of the structure, when $\ell$ is large. Conclusions are made in the last section.

\section{Theoretical foundation for fixed-point iteration method}\label{completeP}

First we give a brief review about the main steps of Pryce's structural analysis method \cite{Pryce2001}.
We consider a DAE system $\textbf{f}=(f_1,f_2,\ldots,f_n)=\textbf{0}$ in $n$ dependent variables $x_{j}=x_{j}(t)$
with $t$ a scalar independent variable, of the form
\begin{equation}
f_i=f_{i}(t,\text{the} \ x_j\ \text{and\ derivatives\ of\ them}), \ i=1,2,\ldots,n.
\end{equation}

Step 1. Form the $n \times n$ signature matrix $\Sigma =
(\sigma_{ij})$ of the DAE, where
\begin{center}
${\sigma_{ij}}=\begin{cases} \text{highest\ differential\ order\ of}\ x_j \ \text{in\ equation}\ f_i, \ \text{if} \ x_j\ \text{appears \ in} \ f_i, \\
 -\infty, \ \text{otherwise}.
\end{cases}$
\end{center}

Step 2. Solve an assignment problem (AP) to find a highest value
transversal (HVT) T, which is a subset of sparsity pattern $S$ with $n$ finite
entries and describes just one element in each row and each column, such that
$\sum\sigma_{ij}$ is maximized and finite. The sparsity pattern $S$
of $\Sigma$ is  defined as:
\begin{equation}
  S=\text{sparse}(\Sigma)=\{(i,j):\sigma_{ij} > -\infty \}.
\end{equation}
This can be formulated as a Linear
Programming Problem (LPP), the Primal is:
\begin{eqnarray}\label{LPP:Primal}
\begin{array}{ll}
  \max\limits_{\xi} & z = \sum\limits_{(i,j)\in S} \sigma_{ij}\xi_{ij},  \label{pequ:a} \\
    \text{s.t.} & \sum\limits_{j:(i,j)\in S} \xi_{ij}=1 \ \text{for\ each} \ i, \label{pequ:b} \\
  & \sum\limits_{i:(i,j)\in S} \xi_{ij}=1 \ \text{for\ each} \ j, \label{pequ:c} \\
   &  \xi_{ij} \geq 0 \ \text{for} \ (i,j) \in S.\label{pequ:d}
\end{array}
\end{eqnarray}

The problem is equivalent to finding a maximum-weight perfect matching in a bipartite graph whose incidence matrix is the signature matrix, and can be solved by Kuhn-Munkres
algorithm\cite{Alexander13} whose time complexity is $O(n^3)$.

Step 3. Determine the offsets of the problem, which are the vectors
$\mathbf{c}= (c_i)_{1\leq i\leq n}, \mathbf{d}=(d_j)_{1\leq j\leq
n}$, the smallest such that $d_j - c_i \geq \sigma_{ij}$, for all
$1\leq i\leq n, 1\leq j\leq n$, and $d_j - c_i = \sigma_{ij}$ when $(i,j)\in T$.

This problem can be formulated as the dual of  (\ref{LPP:Primal}) in the variables $\mathbf{c} = (c_1, c_2, \ldots, c_n$) and $
\mathbf{d} = (d_1, d_2, \ldots, d_n)$. The Dual is defined as follows:
\begin{eqnarray}\label{LPP:Dual}
\begin{array}{ll}
  \min\limits_{\textbf{c},\textbf{d}} & z = \sum\limits_j d_j - \sum\limits_i c_i, \label{equ:a}   \\
  \rm{s.t.} & d_j - c_i \geq \sigma_{ij} \ \text{for\ all} \ (i, j), \\
  & c_i \geq 0 \ \text{for\ all}\ i.
\end{array}
\end{eqnarray}

Step 4. Form the system Jacobian matrix $\mathbf{J}$, given by
\begin{center}
$\mathbf{J}_{ij}=\begin{cases}
 \frac{\partial f_i}{\partial ((d_j-c_i)\text{th\ derivative\ of}\ x_j)}, \ \ \text{if\ this\ derivative\ is\ present\ in}\ f_i,  \\
0, \ \ \ \ \ \ \ \ \ \ \ \ \ \ \ \ \ \ \ \  \ \ \ \ \ \ \ \ \ \ \ \ \
\  \text{otherwise}.
\end{cases}$
\end{center}

Step5. Choose a consistent point. If $\mathbf{J}$ is non-singular at that point, then the solution can be computed with Taylor series or numerical homotopy continuation techniques in a neighborhood of that point.
And using the smallest offsets $\mathbf{c}, \mathbf{d}$ of Problem \ref{LPP:Dual}, the structural index is then defined as:
\begin{equation}\label{SIndex}
\nu=\max_{i}c_i+\begin{cases}
 0, \ \text{for\ all}\ d_j>0, \\
1, \ \text{for\ some}\ d_j=0.
\end{cases}
\end{equation}

In order to determine the smallest offsets of DAEs using fixed-point iteration algorithm \cite{Pryce2001}, we introduce some necessary definitions, firstly.
Define a natural semi-ordering of vectors in $\mathbb{R}^n$, for $\forall \ \textbf{a}, \textbf{b}$, $\textbf{a}\prec\textbf{b}$ if $a_i \leq b_i$ for each $i$, smallest of offsets is in the sense of ordering $\prec$.
Given $\Sigma$ of DAEs systems and a corresponding transversal $T$,
for $\forall \ \textbf{c}=(c_i) (\in \mathbb{R}^n)$, we define a mapping
$$\mathcal{D}(\textbf{c})=(d_j), \mbox{ where } d_j=\max_i (\sigma_{ij}+c_i), $$
and for $\forall \ \textbf{d}=(d_j) (\in \mathbb{R}^n)$, we define a mapping
$$\mathcal{C}_{T}(\textbf{d})=(c^*_i), \mbox{ where } c^*_i=d_j-\sigma_{i,j},\; (i,j) \in T.$$
Furthermore, we define the composition mapping $\phi_{T}(\textbf{c})=\mathcal{C}_{T}(\mathcal{D}(\textbf{c}))$ from $\mathbb{R}^n$ to $\mathbb{R}^n$. Then we obtain fixed-point iteration algorithm below.
\begin{algorithm}[h]        
\caption{Fixed-point iteration algorithm}             
\label{alg:FIFDV}                  
\begin{algorithmic}[1]
\REQUIRE ~~\\       
$\Sigma$ is the signature matrix of DAEs
\ENSURE ~~\\          
$\textbf{c}$ and $\textbf{d}$
\STATE  Set $\textbf{c}' \leftarrow \textbf{0}$
\STATE  Set $\textbf{d} \leftarrow \mathcal{D}(\textbf{c}') $
\STATE  Set T is the HVT of $\Sigma$ computed by Kuhn-Munkres algorithm
\STATE  Set $\textbf{c} \leftarrow \mathcal{C}_{T}(\textbf{d}) $
\WHILE{$\textbf{c} \neq \textbf{c}'$ }
   \STATE Set $\textbf{c}' \leftarrow \textbf{c}$
   \STATE Set $\textbf{d} \leftarrow \mathcal{D}(\textbf{c}) $
   \STATE Set $\textbf{c} \leftarrow \mathcal{C}_{T}(\textbf{d}) $
\ENDWHILE
\RETURN{$\textbf{c}$,$\textbf{d}$}
\end{algorithmic}
\end{algorithm}

In order to give a novel and refined proof of the existence and uniqueness of smallest offsets of Problem \ref{LPP:Dual} and the termination of Algorithm \ref{alg:FIFDV},
we introduce some definitions and lemmas as follows.
\begin{lemma}
\emph{(\cite{Pryce2001})}.
\begin{description}
  \item[(i)] If the Primal has a feasible solution, then it has a basic feasible
  solution (BFS). At such a solution the $\xi_{ij}$ are 1 on some transversal T
  and $0$ elsewhere. The corresponding objective function value is
  $\Sigma_{(i,j)\in T}\sigma_{ij}$, denoted by $||T||$. The optimum is achieved
  at a BFS.
  \item[(ii)] The following results are equivalent.
  \begin{description}
    \item[(a)] The AP is regular.
    \item[(b)] The Primal has a feasible solution.
    \item[(c)] Primal and Dual have a common, finite optimal value given by
    \begin{equation}
      z=\sum d_j-\sum c_i=||T||.
    \end{equation}
  \end{description}
  \item[(iii)(Principle of complementary slackness)] Given a Primal BFS, i.e., a
  transversal T, and a Dual feasible solution \textbf{c},\textbf{d}, the
  following are equivalent:
  \begin{description}\label{Slack}
    \item[(a)] T is an HVT, and \textbf{c} and \textbf{d} are optimal for the Dual.
    \item[(b)] $d_j-c_i=\sigma_{ij}$, for each $(i,j)\in T$.
  \end{description}
\end{description}
\end{lemma}

\begin{lemma}\label{OptFix}\emph{(\cite{Pryce2001})}.
Assume that T is HVT, then $\textbf{c}$ is optimal of Problem \ref{LPP:Dual}
if and only if $\textbf{c}$ is the non-negative fix-points of $\phi_T$,
that is, $\phi_T(\textbf{c})=\textbf{c}$.
\end{lemma}

\begin{definition}
For a given $\Sigma$ and a corresponding T,
the vector set VC is defined as:
\begin{equation}
VC(T)=\{ \textbf{c} \in \mathbb{R}^n | \phi_T(\textbf{c})=\textbf{c}
,\textbf{c} \succ \textbf{0} \}.
\end{equation}
\end{definition}

If T is HVT, the part optimal solution set of Problem \ref{LPP:Dual} is $VC(T)$ by Lemma \ref{OptFix}.
Furthermore, we have:
\begin{lemma}\label{Cindepen}
For a given $\Sigma$ matrix, the optimal-dual set $VC$
is independent of the choice of HVT, that is, $VC(T)=VC(T')$
for any two HVT T and $T'$.
\end{lemma}
\proof
For any $\textbf{c} \in VC(T)$, define  $\textbf{d}= \mathcal{D}(\textbf{c})$.
Then $\textbf{c}$ and $\textbf{d}$ are optimal-dual by Lemma \ref{OptFix}.
Note that $T'$ is HVT. By Lemma \ref{Slack}(iii), obtain
\begin{equation}\label{Cindepen1}
  d_j-c_i=\sigma_{ij},
\end{equation}
for each $(i,j) \in T'$. According to (\ref{Cindepen1}), we have
\begin{equation}\label{}
 \textbf{c}=
\mathcal{C}_{T'}(\textbf{d})=\mathcal{C}_{T'}(\mathcal{D}(\textbf{c}))=
\phi_{T'}(\textbf{c}).
\end{equation}
Therefore, $\textbf{c}$ is also the
non-negative fix-points of $\phi_{T'}$ , that is, $\textbf{c} \in VC(T')$.
Conversely, we can easily prove $VC(T') \subseteq VC(T)$ with the similar
principle above.
\foorp

Now,using the above results we can prove the existence and uniqueness of
the smallest optimal solution for the Dual problem.
\begin{lemma}\label{Smallest Dual}
Assume that the $\Sigma$ matrix of given DAE systems in Problem \ref{LPP:Dual}
contains a transversal T at least, then there exists a unique smallest
dual-optimal pair $\textbf{c}^*$ and $\textbf{d}^*$.
\end{lemma}
\proof
The $\Sigma$ matrix contains a transversal T at least, then there must
exist a HVT T from the finiteness of transversal.
From the Lemma \ref{Cindepen}, assume T is any HVT.
According to the primal-dual principle,
dual-optimal pair $\textbf{c}$ and $\textbf{d}$ must exist, that is,
$VC=VC(T)$ is a non-empty set by Lemma \ref{OptFix}. Moreover,
It is easy to know that for any non-negative vector $\Theta=(\theta,\theta,\ldots,\theta)$,
$\textbf{c}+\Theta$ and $\textbf{d}+\Theta$ is also dual-optimal.
Then $VC$ is a infinite set. Define
\begin{equation}\label{SmallestD1}
\begin{array}{ccc}
  VC_1 = \{||\textbf{c}||_1:\textbf{c} \in VC \}& \text{and}
  & \alpha=inf\{VC_1\}(\geq 0).
\end{array}
\end{equation}

(The existence of `smallest' dual-optimal in ordering `$\leq$')
In fact, the coefficients of all the constraint equations in Problem \ref{LPP:Dual}
are $1$ or $-1$,  and each $\sigma_{ij} \in \Sigma$ is integer. Thus, all the vertices of the feasible region in Linear programming are
integer. Then the `smallest' dual-optimal also are integer.
Therefore, there exists $\textbf{c}^* \in VC$ such that
$\alpha=||\textbf{c}^*||_1$, that is,
$\textbf{c}^*$ and $\textbf{d}^*=\mathcal{D}(\textbf{c}^*)$
is the `smallest' dual-optimal in the sense of ordering `$\leq$'.

(The uniqueness of `smallest' dual-optimal in ordering `$\leq$')
Assume that there are two different `smallest' dual-optimal pair $\textbf{c}^*,\textbf{d}^*$
and $\textbf{c}^o,\textbf{d}^o$ such that
\begin{equation}\label{SmallestD2}
  \alpha=||\textbf{c}^*||_1=||\textbf{c}^o||_1.
\end{equation}
There must exist $i_0 \in \{1,2,\ldots,n\}$ such that
${c^*}_{i_0} \neq {c^o}_{i_0}$. According to the following rules,
construct vector pair $\textbf{c}^{*o}$ and $\textbf{d}^{*o}$.
For given HVT T and each $i$,
if ${c^*}_{i} \geq {c^o}_{i}$, define $ {c^{*o}}_{i}={c^o}_{i}$ and
$ {d^{*o}}_{j}={d^o}_{j} $ such that $(i,j) \in T $;
otherwise, define $ {c^{*o}}_{i}={c^*}_{i}$ and
$ {d^{*o}}_{j}={d^*}_{j} $.

Firstly, it is verified that
$\textbf{c}^{*o}$ and $\textbf{d}^{*o}$ are the Dual feasible solution.
By the definition of $\textbf{c}^{*o}$ and $\textbf{d}^{*o}$,
for each $j \in \{1,2,\ldots,n\}$, obtain
\begin{equation}\label{SmallestD3}
  d^{*o}_{j}=d^{*}_{j} \geq c^{*}_{i}+\sigma_{ij}
\geq c^{*o}_{i}+\sigma_{ij}, i=1,2,\ldots,n;
\end{equation}
and
\begin{equation}\label{SmallestD4}
  d^{*o}_{j}=d^{o}_{j} \geq c^{o}_{i}+\sigma_{ij} \geq c^{*o}_{i}+\sigma_{ij}, i=1,2,\ldots,n.
\end{equation}
Then together with (\ref{SmallestD3}) and (\ref{SmallestD4}),
for each $j \in \{1,2,\ldots,n\}$, get
\begin{equation}\label{SmallestD5}
  d^{*o}_{j} \geq c^{*o}_{i}+\sigma_{ij}, i=1,2,\ldots,n.
\end{equation}
That is $\textbf{c}^{*o}$ and $\textbf{d}^{*o}$ are
the Dual feasible solution.

Furthermore, note that T is HVT,
$\textbf{c}^*,\textbf{d}^*$ and $\textbf{c}^o,\textbf{d}^o$ are
dual-optimal pair.
By Lemma $\ref{Slack}(iii)$, obtain
\begin{eqnarray}\label{SmallestD6}
\begin{array}{cc}
  d^{*}_{j}-c^{*}_{i} = \sigma_{ij},
  & d^{o}_{j}-c^{o}_{i} = \sigma_{ij},
\end{array}
\end{eqnarray}
for each $(i,j) \in T $.
By (\ref{SmallestD6}), we have
\begin{equation}\label{SmallestD7}
  d^{*o}_{j}-c^{*o}_{i}=\sigma_{ij},
\end{equation}
for each $(i,j) \in T $.
Combining (\ref{SmallestD5}) and (\ref{SmallestD7}),
it is indicated that $\textbf{c}^{*o}$ and $\textbf{d}^{*o}$
are also the dual-optimal by Lemma $\ref{Slack}(iii)$,
that is, $\textbf{c}^{*o} \in VC$. But
$||\textbf{c}^{*o}||_1 < ||\textbf{c}^{*}||_1=\alpha$ or
$||\textbf{c}^{*o}||_1 < ||\textbf{c}^{o}||_1=\alpha$,
which is in conflict with (\ref{SmallestD1}).
Therefore, the smallest dual-optimal is unique in the sense of ordering `$\leq$'.

(Smallest in ordering $\prec$ ) Set $\textbf{c}^*$ and $\textbf{d}^*$
are the `smallest' dual-optimal, $\textbf{c}$ and $\textbf{d}$ are any
dual-optimal, then obtain $\alpha=||\textbf{c}^{*}||_1 \leq ||\textbf{c}||_1$.
Assume that there exists $i_0 \in \{1,2,\ldots,n\}$ such that
${c^*}_{i_0} > {c}_{i_0}$. We can construct the new dual-optimal
$\textbf{c}^{o*}(\in VC)$ and $\textbf{d}^{o*}$ such that
$||\textbf{c}^{o*}||_1 <||\textbf{c}^*||_1=\alpha$ by the method described above,
which is also in conflict with (\ref{SmallestD1}).
So obtain $\textbf{c}^* \prec \textbf{c} $, and
$\textbf{d}^*=\mathcal{D}(\textbf{c}^*) \prec \mathcal{D}(\textbf{c})=\textbf{d}$.
Therefore, $\textbf{c}^*$ and $\textbf{d}^*$ is the unique smallest dual-optimal pair.
\foorp

According to the above lemmas, we can prove the termination of
fixed-point iteration algorithm and analyze its complexity.
\begin{lemma}\label{alg:FIFDV_Proof}
The fixed-point iteration algorithm can find
the unique smallest dual-optimal pair $\textbf{c}^*$ and $\textbf{d}^*$ of Problem
\ref{LPP:Dual} by at most $||\textbf{c}^*||_1+1$ iterations.
\end{lemma}
\proof
Set $\textbf{c}^{(1)}=\phi(\textbf{0})(\succ \textbf{0})$,
$\textbf{c}^{(k)}=\phi(\textbf{c}^{(k-1)})=\phi^{k}(\textbf{0})$
for $k \in \mathds{N}^+$. It is verified that $\phi(=\phi_T)$ is monotone operator from
$\mathbb{R}^n$ to $\mathbb{R}^n$, that is, if $\textbf{c} \prec \textbf{c}'$,
$\phi(\textbf{c}) \prec \phi(\textbf{c}')$. So $\{ \textbf{c}^{(k)} \}$ is
a increasing sequence in ``$\prec$" sense, and $\{ ||\textbf{c}^{(k)}||_1 \}$
is also a increasing sequence. Note that T is HVT, then exist
the unique smallest dual-optimal pair $\textbf{c}^*(\succ \textbf{0})$
and $\textbf{d}^*$ by Lemma $\ref{Smallest Dual}$.
According to the monotonicity of $\phi$, obtain
\begin{equation}\label{algProof1}
  \textbf{c}^{(k)} \prec \textbf{c}^{*}, \text{for}  \ k \in \mathds{N}^+,
\end{equation}
and then
$||\textbf{c}^{(k)}||_1 \leq ||\textbf{c}^{*}||_1$. It is indicated
that $\{ ||\textbf{c}^{(k)}||_1 \}$ is bounded.
Based on bounded monotonic principle,
exists $\beta$ such that
\begin{equation}\label{algProof2}
  ||\textbf{c}^{(k)}||_1 \rightarrow \beta (\leq ||\textbf{c}^{*}||_1), k \rightarrow \infty.
\end{equation}

Assume $\{ ||\textbf{c}^{(k)}||_1 \}$ is strictly increasing sequence.
For all $c^{(k)}_i$ are integer, we have
$||\textbf{c}^{(k)}||_1 \geq k$. So $\beta=+\infty$, which is in conflict with (\ref{algProof2}).
It is shown that exists $m \in \mathds{N}^+$
such that $||\textbf{c}^{(m)}||_1=||\textbf{c}^{(m+1)}||_1$. Moreover, obtain
$\textbf{c}^{(m)}=\textbf{c}^{(m+1)}$ and $\textbf{c}^{(l)}=\textbf{c}^{(m)}$ for
$\forall \ l>m$. Set $\textbf{c}^o = \textbf{c}^{(m)}$, from (\ref{algProof1}), we have
$\textbf{c}^o \prec \textbf{c}^*$. In addition, $\textbf{c}^*$ is the smallest
dual-optimal and $\textbf{c}^o$ is the non-negative fix-point,
then we get $\textbf{c}^* \prec \textbf{c}^o$ by Lemma \ref{OptFix} and \ref{Smallest Dual}.
Therefore,
$\textbf{c}^* = \textbf{c}^o$ and $\textbf{d}^* = \mathcal{D}(\textbf{c}^*)$.
Note that the sequence $||\textbf{c}^{(k)}||_1$ increases by 1 on each
iteration at most. So the Algorithm \ref{alg:FIFDV} finds
the smallest dual-optimal pair $\textbf{c}^*$ and $\textbf{d}^*$ by at most $||\textbf{c}^*||_1+1$ iterations.
\foorp

Now, we are able to give the main theorem below.
\begin{theorem}\label{PryceTime}
If the $\Sigma$ matrix of given DAE systems in Problem \ref{LPP:Dual}
contains some transversal $T$, then
the unique smallest dual-optimal pair $\textbf{c}^*$ and $\textbf{d}^*$
can be found in time $O(n^3+||\textbf{c}^*||_1 \cdot n^2)$ via
fixed-point iteration algorithm.
\end{theorem}
\proof
It can be easily proved by Lemma \ref{Smallest Dual} and \ref{alg:FIFDV_Proof}.
\foorp

\begin{example}\label{ex1}
Consider the application of the Algorithm \ref{alg:FIFDV} to nonlinear
DAE system $\textbf{f}=(f_1,f_2,f_3)=\textbf{0}$ in three dependent variables $x_1(t), x_2(t), x_3(t)$
with known forcing functions $u_{i}(t)(i=1,2,3)$:
\begin{eqnarray*}
\left\{\aligned
f_1 =& \ddot{x_1}+x_3+u_{1}(t) \\
f_2 =& \dot{x_2}+x_3+u_{2}(t) \\
f_3 =& {x_1}^2+{x_2}^2+u_{3}(t)
\endaligned\right..\end{eqnarray*}
\end{example}

The corresponding signature matrix is
\begin{equation*}
  \Sigma=\left(
      \begin{array}{ccc}
        2^* &  & 0 \\
         & 1 & 0^* \\
        0 & 0^* &  \\
      \end{array}
    \right),
\end{equation*}
where we have already marked the HVT with asterisks, and the elements in the blanks of $\Sigma$ are $-\infty$. We give the
main process of Algorithm \ref{alg:FIFDV} below,
\begin{equation*}
\begin{array}{ll}
\Sigma \Rightarrow \begin{blockarray}{ccccc}
~ & x_{1} &  x_{2} &  x_{3}  &  {\textbf{c}'}^{(0)}\\
\begin{block}{c[ccc]c}
f_{1} & 2^* & ~ & 0 & 0 \\
f_{2} & ~ & 1 & 0^* & 0 \\
f_{3} & 0 & 0^* & ~ & 0 \\
\end{block}
\end{blockarray}
& \Rightarrow\begin{blockarray}{cccccc}
 ~ & x_{1} & x_{2} & x_{3} & {\textbf{c}'}^{(0)} & \textbf{c}^{(0)} \\
\begin{block}{c[ccc]cc}
f_{1} & 2^* & ~  & 0 & 0 & 0\\
f_{2} & ~ & 1 & 0^* & 0 & 0\\
f_{3} & 0 & 0^* & ~  & 0 & 1 \\
\end{block}
{\textbf{d}}^{(0)} & 2 & 1 & 0 & \\
\end{blockarray} \\
\end{array}
\end{equation*}
where ${\textbf{c}'}^{(i)}$, $\textbf{c}^{(i)}$, and ${\textbf{d}}^{(i)}$ mean the $i$th
iteration for $\textbf{c}'$, $\textbf{c}$ and $\textbf{d}$, respectively.
Therefore, obtain the smallest offsets $\textbf{c}=(0,0,1)$ and
$\textbf{d}=(2,1,0)$ of the DAE systems.

\section{Block fixed-point iteration method}\label{SectBFIRM}

When dealing with DAE systems of large dimensions, an important manipulation is the block triangularization of the system \cite{Maff96}, which allows to decompose the overall system into subsystems which can be solved in sequence.
Similarly, considering the index reduction for large-scale systems, it is necessary to compute
the block triangular forms (BTF) of the $\Sigma$ matrix by permuting its rows and columns \cite{Pothen90,Duff10}.

In this section, assume the given DAE systems are structurally nonsingular, meaning that the $\Sigma$ matrix of the systems exists a transversal, then obtain the BTF of
the $\Sigma$ matrix below,
\begin{equation}\label{MBTF}
M=
\begin{pmatrix}
  M_{11} & M_{1,2} & \cdots & M_{1,\ell} \\
         & M_{2,2} & \cdots & M_{2,\ell} \\
         &  & \ddots & \vdots \\
         &  &  & M_{\ell,\ell}
\end{pmatrix},
\end{equation}
where the elements in the blanks of $M$ are $-\infty$ ,
and diagonal matrix $M_{i,i}$ is square and irreducible, for $i=1,2,\ldots,\ell$ \cite{Pryce2012}.

The main idea of block fixed-point iteration method for $\Sigma$ matrix
with BTF is to use the fixed-point iteration method with parameter mentioned below to process each diagonal matrix in block upper triangulated signature matrix from top to bottom in
sequence. We give the fixed-point iteration method with parameter, firstly.

\subsection{Fixed-point iteration method with parameter}
The Dual Problem \ref{LPP:Dual} with $n$ dimension nonnegative parameter vector $\textbf{p}$ is defined as follows:
\begin{definition}
The Dual Problem with nonnegative parameter $\textbf{p}$ is defined via:
\begin{eqnarray}\label{Dual_md}
\begin{array}{ll}
  \min\limits_{\textbf{c},\textbf{d}} & z = \sum\limits_j d_j - \sum\limits_i c_i,    \\
  \text{s.t.} & d_j - c_i \geq \sigma_{ij} \ \text{for\ all}\ (i, j),  \\
  &d_j \geq p_j, \ \text{for\ all}\ j, \\
  & c_i \geq 0 \ \text{for\ all}\ i. \\
\end{array}
\end{eqnarray}
\end{definition}

For any given nonnegative parameter $\textbf{p}$, obtain the fixed-point iteration
algorithm with parameter (PFPIA) below just by modified
the fixed-point iteration algorithm.
\begin{algorithm}[!h]        
\caption{Fixed-point iteration algorithm with parameter (PFPIA)}       
\label{alg:FIFDV_MD}         
\begin{algorithmic}[1]
\REQUIRE ~~\\       
$\Sigma$ is signature matrix for DAE systems

$\textbf{p}$ is nonnegative parameter vector
\ENSURE ~~\\          
$\textbf{c}$ and $\textbf{d}$
\STATE  Set $\textbf{c}' \leftarrow \textbf{0}$
\STATE  Set T is HVT of $\Sigma$ by Kuhn-Munkres algorithm
\STATE  Set $\textbf{c}' \leftarrow \mathcal{C}_{T}(\textbf{p})$
\STATE  Set $\textbf{c}' \leftarrow \max \{ \textbf{c}',\textbf{0} \}$
\STATE  Set $\textbf{d} \leftarrow \mathcal{D}(\textbf{c}') $
\STATE  Set $\textbf{c} \leftarrow \mathcal{C}_{T}(\textbf{d}) $
\WHILE{$\textbf{c} \neq \textbf{c}'$ }
   \STATE Set $\textbf{c}' \leftarrow \textbf{c}$
   \STATE Set $\textbf{d} \leftarrow \mathcal{D}(\textbf{c}) $
   \STATE Set $\textbf{c} \leftarrow \mathcal{C}_{T}(\textbf{d}) $
\ENDWHILE
\RETURN{$\textbf{c}$,$\textbf{d}$}
\end{algorithmic}
\end{algorithm}

\begin{lemma}\label{Smallest Dual_MD}
Let $\textbf{p}$ is any nonnegative parameter.
Assume that the $\Sigma$ matrix in the Problem \ref{Dual_md} contains a transversal
T, then there exists a unique smallest dual-optimal pair
$\textbf{c}^*$ and $\textbf{d}^*$ such that
\begin{equation*}
\textbf{c}^*=\min \{ \textbf{c}|\mathcal{D}(\textbf{c})\succ \textbf{p}
,\textbf{c} \in VC\}.
\end{equation*}
Moreover, if the used transversal T is HVT, the Algorithm
\ref{alg:FIFDV_MD} finds
the unique smallest dual-optimal pair $\textbf{c}^*$ and $\textbf{d}^*$ by at most $||\textbf{c}^*||_1 - ||\max\{\mathcal{C}_T(\textbf{p}),\textbf{0}\}||_1+1$ iterations.
\end{lemma}
\proof
Just modify the proof of Lemma \ref{Smallest Dual} and \ref{alg:FIFDV_Proof}
properly.
\foorp

\begin{remark}\label{DuMdtoDu}
If $p_j \leq \max_{i}\sigma_{ij}$, for $j=1,2,\ldots,n$, and
$d_j \geq \max_{i}\sigma_{ij}$ derived from Problem (\ref{Dual_md}),
for $j=1,2,\ldots,n$. So obtain $d_j \geq p_j$, for $j=1,2,\ldots,n$,
that is, the constraint condition $d_j \geq p_j$ in Problem (\ref{Dual_md}) can be
deleted. Therefore, the Problem (\ref{Dual_md}) turns into Problem (\ref{LPP:Dual}).
\end{remark}

\begin{example}\label{ex2}
Consider the application of the Algorithm \ref{alg:FIFDV_MD} to nonlinear
DAE systems $\textbf{f}=(f_4,f_5,f_6)=\textbf{0}$ in three dependent variables $x_4(t), x_5(t) ,x_6(t)$
with known forcing functions $u_{i}(t)(i=4,5,6)$:
\begin{eqnarray*}
\left\{\aligned
f_4 = &\ddot{x_4}+x_6+u_{4}(t) \\
f_5 = &\dot{x_5}+x_6+u_{5}(t) \\
f_6 = &{x_4}^2+{x_5}^2 \\
\endaligned\right.,\end{eqnarray*}
and the any given parameter is $\textbf{p}=(0,0,2)$.
\end{example}
We give below the main process of Algorithm \ref{alg:FIFDV_MD},
\begin{equation*}
\begin{array}{ll}
\Sigma \Rightarrow\begin{blockarray}{ccccc}
\textbf{p} & 0 & 0 & 2 & \\
~ & x_{4} &  x_{5} &  x_{6}  &  {\textbf{c}'}^{(0)} \\
\begin{block}{c[ccc]c}
f_{4} & 2^* & ~ & 0 & 0 \\
f_{5} & ~ & 1 & 0^* & 2 \\
f_{6} & 0 & 0^* & ~ & 0 \\
\end{block}
\end{blockarray} &
\Rightarrow\begin{blockarray}{cccccc}
\textbf{p} & 0 & 0 & 2 & \\
 ~ & x_{4} & x_{5} & x_{6} & {\textbf{c}'}^{(0)} &  {\textbf{c}}^{(0)}\\
\begin{block}{c[ccc]cc}
f_{4} & 2^* & ~  & 0 & 0 & 0\\
f_{5} & ~ & 1 & 0^* & 2 & 2\\
f_{6} & 0 & 0^* & ~  & 0 & 3 \\
\end{block}
{\textbf{d}}^{(0)} & 2 & 3 & 2 & \\
\end{blockarray} \\
\Rightarrow\begin{blockarray}{ccccccc}
\textbf{p} & 0 & 0 & 2 & &\\
 ~ & x_{4} & x_{5} & x_{6} & {\textbf{c}'}^{(0)} &  {\textbf{c}}^{(0)}&{\textbf{c}}^{(1)} \\
\begin{block}{c[ccc]ccc}
f_{4} & 2^* & ~  & 0 & 0 & 0 &1 \\
f_{5} & ~ & 1 & 0^* & 2  & 2 &2 \\
f_{6} & 0 & 0^* & ~  & 0 & 3 &3 \\
\end{block}
{\textbf{d}}^{(0)}    & 2 & 3 & 2 & & \\
{\textbf{d}}^{(1)} & 3 & 3 & 2 & & \\
\end{blockarray} &
\Rightarrow\begin{blockarray}{cccccccc}
\textbf{p} & 0 & 0 & 2 & & &\\
 ~ & x_{4} & x_{5} & x_{6} &{\textbf{c}'}^{(0)} &  {\textbf{c}}^{(0)}&{\textbf{c}}^{(1)}&{\textbf{c}}^{(2)} \\
\begin{block}{c[ccc]cccc}
f_{4} & 2^* & ~  & 0 & 0 & 0 &1 &1 \\
f_{5} & ~ & 1 & 0^* &  2  & 2 &2 &2 \\
f_{6} & 0 & 0^* & ~  & 0 & 3 &3 &3\\
\end{block}
{\textbf{d}}^{(0)}& 2 & 3 & 2 & & \\
{\textbf{d}}^{(1)}& 3 & 3 & 2 & & \\
{\textbf{d}}^{(2)}& 3 & 3 & 2 & & \\
\end{blockarray},
\end{array}
\end{equation*}
where ${\textbf{c}'}^{(i)}$, $\textbf{c}^{(i)}$ and ${\textbf{d}}^{(i)}$ mean the $i$th
iteration for $\textbf{c}'$, $\textbf{c}$ and $\textbf{d}$, respectively.
Then obtain the smallest offsets $\textbf{c}=(1,2,3)$ and
$\textbf{d}=(3,3,2)$ for the DAEs.

\subsection{Block fixed-point iteration method}
The given DAE systems are structurally nonsingular,
obtain the $\Sigma$ matrix $M$ of Problem \ref{LPP:Dual}
with block triangular form (\ref{MBTF}),
and $\sum\limits_{i=1}^{\ell}n_i=n$, where $n_i$ is the order of
$M_{ii},i=1,2,\ldots,\ell$.
In order to find the unique smallest dual-optimal,
we give some necessary symbols and definitions as follows.

Let the parameter vector
is $\textbf{p}=(\textbf{p}_1,\textbf{p}_2,\ldots,\textbf{p}_\ell)$
with $\ell$ sections, the dual-optimal are
$\textbf{c}=(\textbf{c}_1,\textbf{c}_2,\ldots,\textbf{c}_\ell)$ and
$\textbf{d}=(\textbf{d}_1,\textbf{d}_2,\ldots,\textbf{d}_\ell)$,
where the dimension of $\textbf{p}_i$, $\textbf{c}_i$ and $\textbf{d}_i$
are $n_i$ for $i=1,2,\ldots,\ell$.
For any $n \times r$ order matrix $B$ and $B'$ ,$n$ order vector $\textbf{q}$,
$\bar{\textbf{q}}$ and $\hat{\textbf{q}}$,
$r$ order vector $\textbf{w}$,
the mapping $B'=\textbf{RowAdd}(B,\textbf{q})$ is defined as $B'_{i,j}=B_{i,j}+q_i$,
for $i=1,2,\ldots,n, j=1,2,\ldots,r$; the mapping $\textbf{w}=\textbf{ColMax}(B)$ is
defined via $w_j=\max\limits_{i\in\{1,2,\ldots,n\}}{B_{i,j}}$, for $j=1,2,\ldots,r$;
the mapping $\textbf{q}=\textbf{eMax}(\bar{\textbf{q}},\hat{\textbf{q}})$ is defined
as $q_i=max(\bar{q}_i,\hat{q}_i)$, $i=1,2,\ldots,n$.
Then we give block fixed-point iteration algorithm below.
\begin{algorithm}[h]        
\caption{Block fixed-point iteration algorithm}      
\label{alg:BFIRA}       
\begin{algorithmic}[1]
\REQUIRE ~~\\       
M is $\Sigma$ matrix of given DAE systems with BTF (\ref{MBTF})
\ENSURE ~~\\          
$\textbf{c}=(\textbf{c}_1,\textbf{c}_2,\ldots,\textbf{c}_\ell)$
and $\textbf{d}=(\textbf{d}_1,\textbf{d}_2,\ldots,\textbf{d}_\ell)$
\STATE  Set $\textbf{p}=(\textbf{p}_1,\textbf{p}_2,\ldots,\textbf{p}_\ell)$
,$\textbf{p}_j=\textbf{0}$,$j=1,2,\ldots,\ell$.
\STATE Set $\textbf{c}=(\textbf{c}_1,\textbf{c}_2,\ldots,\textbf{c}_\ell)$,
$\textbf{c}_j=\textbf{0}$,$j=1,2,\ldots,\ell$
\STATE Set $\textbf{d}=(\textbf{d}_1,\textbf{d}_2,\ldots,\textbf{d}_\ell)$
,$\textbf{d}_j=\textbf{0}$,$j=1,2,\ldots,\ell$.
\STATE Get $(\textbf{c}_1,\textbf{d}_1)=\text{PFPIA}(M_{11},\textbf{p}_1)$.
\FOR{$i \ \text{from} \ 2 \ \text{to} \ \ell $}
    \STATE Update: $[M_{i-1,i},\ldots,M_{i-1,\ell}] \leftarrow
    \textbf{RowAdd}([M_{i-1,i},\ldots,M_{i-1,\ell}],\textbf{c}_{i-1})$.
    \STATE Update:$\textbf{p}_i \leftarrow
    \textbf{ColMax}(
    \begin{pmatrix} M_{1,i} \\  \vdots \\ M_{i-1,i} \end{pmatrix})$.
    \STATE  Update: $\textbf{p}_i \leftarrow \textbf{eMax}(\textbf{p}_i,\textbf{0})$.
    \STATE Get $(\textbf{c}_i,\textbf{d}_i)=\text{PFPIA}(M_{ii},\textbf{p}_i)$.
\ENDFOR
\RETURN $\textbf{c}$,$\textbf{d}$
\end{algorithmic}
\end{algorithm}

In order to obtain a complete theoretical analysis of
block fixed-point iteration method,
we give some necessary lemmas, firstly.
\begin{lemma}\label{UnionTrans}
\emph{(\cite{Pryce2012})}.
For given $\Sigma$ matrix $M$ with BTF (\ref{MBTF}),
if $T_i$ is HVT of $M_{ii}$, then
$T=\bigcup\limits_{i=1}^{\ell}T_i$ is HVT of $M$.
\end{lemma}

\begin{lemma}\label{PryceEBlock}
Assume that the Problem \ref{LPP:Dual} of DAE systems is structurally nonsingular,
then fixed-point iteration algorithm gives the same smallest
dual-optimal pair with block fixed-point iteration algorithm.
\end{lemma}
\proof
Without loss of generality, assume that the $\Sigma$ matrix $M$ of Problem \ref{LPP:Dual}
is block triangular forms (\ref{MBTF}). Set
$\textbf{p}=(\textbf{p}_1,\textbf{p}_2,\ldots,\textbf{p}_\ell)$
with $\ell$ sections is the parameter vector;
$\textbf{c}^o=(\textbf{c}^{o}_1,\textbf{c}^{o}_2,\ldots,\textbf{c}^{o}_\ell)$ and
$\textbf{d}^o=(\textbf{d}^{o}_1,\textbf{d}^{o}_2,\ldots,\textbf{d}^{o}_\ell)$
are the smallest dual-optimal found by block fixed-point iteration algorithm;
$\textbf{c}^*=(\textbf{c}^{*}_1,\textbf{c}^{*}_2,\ldots,\textbf{c}^{*}_\ell)$ and
$\textbf{d}^*=(\textbf{d}^{*}_1,\textbf{d}^{*}_2,\ldots,\textbf{d}^{*}_\ell)$
are the smallest dual-optimal found by fixed-point iteration algorithm.

For $\ell$ is integer, we prove the lemma by mathematical induction. Considering about
$\ell=1$, it is easy to know $\textbf{p}_1=\textbf{0}$, so then
$\textbf{c}^{o}=\textbf{c}^{*}$ and $\textbf{d}^{o}=\textbf{d}^{*}$, that is,
the lemma is true.
Assume the lemma is true when $\ell=N-1$, that is,
\begin{equation}\label{BlockP1}
  \begin{array}{cc}
    \textbf{c}^{o}_k=\textbf{c}^{*}_k,
     & \textbf{d}^{o}_k=\textbf{d}^{*}_k,
  \end{array}
\end{equation}
for $k=1,2,\ldots,N-1$.
We now consider about $\ell=N$. From (\ref{BlockP1}), obtain
\begin{equation}\label{BlockP2}
\begin{array}{rl}
\textbf{p}_N &=\textbf{eMax}(\textbf{ColMax}(\textbf{RowAdd}(
\begin{pmatrix} M_{1,N} \\  \vdots \\ M_{N-1,N} \end{pmatrix},
\begin{pmatrix} \textbf{c}^{o}_1 \\  \vdots \\ \textbf{c}^{o}_{N-1}
\end{pmatrix}) ),\textbf{0}) \\
  &=\textbf{eMax}(\textbf{ColMax}(\textbf{RowAdd}(
\begin{pmatrix} M_{1,N} \\  \vdots \\ M_{N-1,N} \end{pmatrix},
\begin{pmatrix} \textbf{c}^{*}_1 \\  \vdots \\ \textbf{c}^{*}_{N-1}
\end{pmatrix})),\textbf{0}).
\end{array}
\end{equation}

It is easily verified that $\textbf{c}^{*}_N$ and $\textbf{d}^{*}_N$ are the dual-optimal
of Problem \ref{Dual_md} with parameter $\textbf{p}_N$ by (\ref{BlockP2}).
From block fixed-point iteration algorithm, note that
\begin{equation}\label{Bc}
\textbf{c}^o_N=\min \{ \textbf{c}|\mathcal{D}(\textbf{c})\succ \textbf{p}_N
,\textbf{c} \in VC\}.
\end{equation}
From $\text{Lemma} \ \ref{Smallest Dual_MD}$ and (\ref{Bc}), we obtain
\begin{equation}\label{BlockP3}
  \begin{array}{cc}
    \textbf{c}^{o}_N \prec \textbf{c}^{*}_N,
     & \textbf{d}^{o}_N=\mathcal{D}(\textbf{c}^{o}_N ) \prec
     \textbf{d}^{*}_N=\mathcal{D}(\textbf{c}^{*}_N ).
  \end{array}
\end{equation}

On the other hand, construct $\textbf{c}^{*o}=(\textbf{c}^{*}_1,\ldots,
\textbf{c}^{*}_{N-1},\textbf{c}^{o}_N)$ and
$\textbf{d}^{*o}=(\textbf{d}^{*}_1,\ldots,
\textbf{d}^{*}_{N-1},\textbf{d}^{o}_N)$.
Note that $\textbf{c}^{o}_N$ and $\textbf{d}^{o}_N$ are the smallest dual-optimal
of Problem \ref{Dual_md} with parameter $\textbf{p}_N$.
Moreover, from (\ref{BlockP1}) and (\ref{BlockP2}),
$\textbf{c}^{*o}$ and $\textbf{d}^{*o}$ are the dual feasible solution
of Problem \ref{LPP:Dual}.
By Lemma \ref{Slack}(iii), obtain
\begin{equation}\label{BlockP4}
  \textbf{d}^{*}_{k_{j_k}}-\textbf{c}^{*}_{k_{i_k}}=\sigma_{i_k,j_k},
\end{equation}
for each $(i_k,j_k) \in T_k$, $k=1,2,\ldots,N-1$, and
\begin{equation}\label{BlockP5}
  \textbf{d}^{o}_{N_{j_N}}-\textbf{c}^{o}_{N_{i_N}}=\sigma_{i_N,j_N}
\end{equation}
for each $(i_N,j_N) \in T_N$.
From (\ref{BlockP4},\ref{BlockP5}) and Lemma \ref{UnionTrans}, we have
\begin{equation}\
  \textbf{d}^{*o}_j-\textbf{c}^{*o}_i=\sigma_{i,j}
\end{equation}
for each $(i,j) \in T=\bigcup\limits_{i=1}^{N}T_i$. That is,
$\textbf{c}^{*o}$ and $\textbf{d}^{*o}$ are dual-optimal of Problem \ref{LPP:Dual}
with $\Sigma=M$ by $\text{Lemma} \ \ref{Slack}\text{(iii)}$.
Note that $\textbf{c}^{*}$ and $\textbf{d}^{*}$ are the smallest dual-optimal
of Problem \ref{LPP:Dual}.
By Lemma \ref{Smallest Dual}, obtain
\begin{equation}\label{BlockP6}
  \textbf{c}^{*} \prec \textbf{c}^{*o}, \textbf{d}^{*} \prec \textbf{d}^{*o}.
\end{equation}
Moreover, by (\ref{BlockP6}), we get
\begin{equation}\label{BlockP7}
  \textbf{c}^{*}_N \prec \textbf{c}^{o}_N,  \textbf{d}^{*}_N \prec \textbf{d}^{o}_N.
\end{equation}
Combining (\ref{BlockP3}) with (\ref{BlockP7}), obtain
\begin{equation}\label{BlockP8}
  \textbf{c}^{o}_N = \textbf{c}^{*}_N, \textbf{d}^{o}_N = \textbf{d}^{*}_N.
\end{equation}
So we have
$\textbf{c}^{*} = \textbf{c}^{o}$
and $\textbf{d}^{*} = \textbf{d}^{o}$.
It is shown that the lemma is true.
\foorp

Now, we are able to obtain the main theorem as follows.
\begin{theorem}\label{BFIRM_Time}
If the $\Sigma$ matrix for given DAE systems with BTF(\ref{MBTF}) is structurally nonsingular,
then the unique smallest dual-optimal pair
$\textbf{c}^*=(\textbf{c}^*_1,\textbf{c}^*_2,\ldots,\textbf{c}^*_\ell)$
and $\textbf{d}^*$ of the DAE can be found in time
$O(\sum\limits_{i=1}^{\ell}{{n_i}^3+||\textbf{c}^*_i||_1 \cdot {n_i}^2})$
by block fixed-point iteration algorithm.
Furthermore, if $n_i=r$ for each $i$, i.e., $n=\ell \cdot r$, then
the time is $O(\ell \cdot r^3+||\textbf{c}^*||_1 \cdot r^2)$.
\end{theorem}
\proof
We can easily prove the theorem by Theorem \ref{PryceTime} and
Lemma \ref{Smallest Dual_MD} and \ref{PryceEBlock}.
\foorp

\begin{example}\label{CounterEx}
Consider the application of the Algorithm \ref{alg:BFIRA} to nonlinear DAE system
$\textbf{f}=(f_1, f_2,\ldots, f_6)=\textbf{0}$ in six dependent variables $x_1(t), x_2(t), \ldots,x_6(t)$ with known forcing functions $u_{i}(t)(i=1,2,\ldots,6)$:
\begin{eqnarray*}
\left\{\aligned
f_1 =& \ddot{x_1}+x_3+u_{1}(t) \\
f_2 =& \dot{x_2}+x_3+u_{2}(t) \\
f_3 =& {x_1}^2+{x_2}^2+\dot{x_6}+u_{3}(t)\\
f_4 =& \ddot{x_4}+x_6 +u_{4}(t)\\
f_5 =& \dot{x_5}+x_6 +u_{5}(t)\\
f_6 =& {x_4}^2+{x_5}^2+u_{6}(t)
\endaligned\right..
\end{eqnarray*}
\end{example}

The corresponding signature matrix is
\begin{equation*}
\Sigma = \begin{blockarray}{ccccccccc}
\textbf{p} & 0 & 0 & 0 & 0 & 0 & 2 & & \\
  ~ & x_{1} &  x_{2}  &  x_{3}  &  x_{4}  &  x_{5}  &  x_{6}  & \tilde{\textbf{c}}^* & \textbf{c}^* \\
\begin{block}{c(ccc|ccc)cc}
f_{1} &  2^* &   & 0 &   &   &    &0 &0  \\
f_{2} &    & 1 & 0^* &   &   &    &0 &0  \\
f_{3} & 0  & 0^* &   &   &   & 1  &1 &1  \\
\cline{2-7}
f_{4} &    &   &   & 2^* &   &0   &0 &1  \\
f_{5} &    &   &   &   & 1 &0^*   &0 &2  \\
f_{6} &    &   &   & 0 & 0^* &   &1 &3  \\
\end{block}
\tilde{\textbf{d}}^* & 2 & 1 & 0 &  2 & 1 & 0 & & \\
        \textbf{d}^* & 2 & 1 & 0 &  3 & 3 & 2 & & \\
\end{blockarray},
\end{equation*}
where we have already marked the $HVT$ with asterisks; $\textbf{p}=(0,0,0,0,0,2)$ is the corresponding parameter vector by block fixed-point iteration algorithm;
$\tilde{\textbf{c}}_i^*=(0,0,1)$ and $\tilde{\textbf{d}}_i^*=(2,1,0)$
are the local smallest offsets for each diagonal signature matrix
$\Sigma_{ii}$, $i=1,2$ via fixed-point iteration algorithm;
$\textbf{c}^*=(0,0,1,1,2,3)$ and $\textbf{d}^*=(2,1,0,3,3,2)$ are the global
smallest offsets for $\Sigma$ matrix directly using fixed-point iteration algorithm.

In the following, the main process of block fixed-point iteration algorithm is shown.
The signature matrix $\Sigma$ above contains two blocks. For $i=1$,
$\textbf{p}_1=(0,0,0)$ is the parameter for the first diagonal block,
get $\textbf{c}_1^*=\tilde{\textbf{c}}_1^*=(0,0,1)$ and
$\textbf{d}_1^*=\tilde{\textbf{d}}_1^*=(2,1,0)$ from Example \ref{ex1}.
For $i=2$, we obtain $\textbf{p}_2=(0,0,2)$ which is the parameter for the second diagonal block, and obtain $\textbf{c}_2^*=(1,2,3)$,
$\textbf{d}_2^*=(3,3,2)$ from Example \ref{ex2}. So
$\textbf{c}^*=(\textbf{c}^*_1,\textbf{c}^*_2)=(0,0,1,1,2,3)$ and
$\textbf{d}^*=(\textbf{d}^*_1,\textbf{d}^*_2)=(2,1,0,3,3,2)$
are the smallest offsets for the DAE system.

\section{Conclusions}\label{Conclusion}

In this paper, we reinforce the theoretical foundation for Pryce's
structural index reduction method of DAE systems, finely prove the existence and uniqueness of the smallest offsets,
and then show the polynomial complexity for finding optimal index reduction for given DAEs.

To solve large scale DAE systems with block structure, we describe a block fixed-point iteration method which can be applied to a sequence of sub-systems rather than the whole system. Accordingly, the time complexity of our method decreases proportionally with the number of the diagonal blocks in the signature matrix.

As pointed in the Campbell-Griepentrog Robot Arm \cite{Griepentrog95} and the special DAE with parameter \cite{Lamour12}, Pryce's structural analysis method fails to find a DAE's true structure because of producing an identically singular Jacobian.
What's more, for a class of simple DAE systems with special $n \times n$ signature matrix $\Sigma$ \cite{Reissig2000}, the actual number of iterations of fixed-point iteration algorithm
(i.e., $O(n)$) is significantly less than $||\textbf{c}^*||_1$.
We believe these situations have appeared rarely in the practical applications.
Compared with other structural index reduction methods,
our  method can address a fairly wide class of large-scale
DAE systems precisely and efficiently.
And the actual performance of block fixed-point iteration algorithm will be discussed
in future work.

\section*{Acknowledgements}
This work is partially supported by China 973 Project (Grant No.
NKBRPC-2011CB302402), National Natural Science Foundation of China (Grant Nos. 11471307, 61402537, 11171053, 91118001), the West Light
Foundation of the Chinese Academy of Sciences from China, the project of Chongqing Science and Technology Commission from China (Grant No. cstc2013jjys40001).

\end{document}